\documentclass[review]{elsarticle}
\usepackage{xcolor}

\usepackage{latexsym,amsmath,amsfonts,amssymb,mathrsfs,bm,appendix,multirow}
\usepackage{float}
\usepackage{booktabs}
\usepackage{amssymb,amsmath,multirow,graphicx}
\usepackage{multicol}    
\usepackage{graphicx}
\usepackage{marvosym}
\usepackage{ifsym}
\usepackage{amssymb,epsfig,multirow}
\usepackage{graphicx}
\usepackage{float}
 \usepackage{color}
\usepackage{enumerate}
\usepackage{xcolor}
\usepackage{amsmath}
\usepackage{ulem}
\usepackage{amsmath}\allowdisplaybreaks
\usepackage{arydshln}
\usepackage{algorithm}
\usepackage{algpseudocode}
\usepackage{amsmath}

\usepackage{lineno,hyperref}

\newtheorem{theorem}{Theorem}[section]

\newtheorem{lemma}[theorem]{Lemma}
\newtheorem{proposition}[theorem]{Proposition}
\newtheorem{definition}[theorem]{Definition}

\newtheorem{remark}[theorem]{Remark}
\modulolinenumbers[5]

\journal{Journal of \LaTeX\ Templates}

\begin{document}

\begin{frontmatter}

\title{A polynomially solvable case of unconstrained ($-$1,1)-quadratic fractional optimization}


\author[label1]{Meijia Yang}
\ead{myang@ustb.edu.cn}
\address[label1]{School of Mathematics and Physics, University of Science and Technology Beijing, Beijing, 100083, P. R. China}
\author[label2]{Yong Xia\corref{cor1}}\cortext[cor1]{Corresponding author.}
\ead{yxia@buaa.edu.cn}
\address[label2]{LMIB of the Ministry of Education, School of Mathematical Sciences, Beihang University, Beijing, 100191, P. R. China}

\begin{abstract}
In this paper, we consider an unconstrained ($-$1,1)-quadratic fractional optimization in the following form: $\min_{x\in\{-1,1\}^n}~(x^TAx+\alpha)/(x^TBx+\beta)$, where $A$ and $B$, given by their nonzero eigenvalues and associated eigenvectors, have ranks not exceeding fixed integers $r_a$ and $r_b$, respectively.  
We show that this problem can be solved in $O(n^{r_a+r_b+1}\log^2 n)$ by the accelerated Newton-Dinkelbach method when the matrices $A$ has nonpositive diagonal entries only, $B$ has nonnegative diagonal entries only. Furthermore, this problem can be solved in $O(n^{r_a+r_b+2}\log^2 n)$ when $A$ has $O(\log(n))$ positive diagonal entries, $B$ has $O(\log(n))$ negative diagonal entries.
\end{abstract}


\end{frontmatter}



\section{Introduction}
Consider the unconstrained ($-$1,1)-quadratic fractional optimization in the following form
\begin{eqnarray*}
{\rm (P)}~~&\min& \frac{x^TAx+\alpha}{x^TBx+\beta} \\
&{\rm s.t.}& x\in \left\{-1,~1\right\}^n,
\end{eqnarray*}
where $A\in\Bbb R^{n\times n}$, $B\in\Bbb R^{n\times n}$ are symmetric matrices. We assume that $x^TBx+\beta>0$ for $\forall x\in\{-1,1\}^n$. {\rm (P)} has a finite optimum objective value, denoted by $\delta^*$. Without loss of generality, we can further assume that $\alpha$ is large enough such that $\delta^*\geq0$. Since otherwise, as we have assumed that $x^TBx+\beta>0$ holds for $\forall x\in\{-1,1\}^n$, we can add $\gamma>0$ to the objective function of {\rm (P)} such that 
\[
0\leq\frac{x^TAx+\alpha}{x^TBx+\beta}+\gamma=\frac{x^TAx+\alpha+\gamma(x^TBx+\beta)}{x^TBx+\beta},
\]
which means that the numerator $x^T(A+\gamma B)x+\alpha+\gamma\beta\geq0$.

{\rm (P)} has many applications in the real world such as image segmentation \cite{YGWD2023}. In \cite{SM2000}, J. Shi and J. Malik proposed a new graph-theoretic criterion for measuring the goodness of an image partition and the model is finally shown to be {\rm (P)} with another constraint that $x^Te=0$. Furthermore, when $A\preceq 0$, $B\succeq 0$, {\rm (P)} is minimizing a concave/covex type fractional problem. Then the problem {\rm (P)} is equivalent to the problem of minimizing $(x^TAx+\alpha)/(x^TBx+\beta)$ over the 
 box constraint $x\in\left[-1,1\right]^n$, which has many applications in image processing \cite{BT2010}.


The linear fractional combinatorial optimization problem {\rm (LFC)} 
has been well studied: 
\begin{eqnarray*}
{\rm (LFC)}~~&\inf& \frac{c^Tx}{d^Tx} \\
&{\rm s.t.}& x\in \mathcal D\subseteq \left\{0,~1\right\}^n,
\end{eqnarray*}
where $c\in\Bbb R^n$, and $d\in\Bbb R^n$ such that $d^Tx>0$ holds for $\forall x\in \mathcal D$. {\rm (LFC)} can be reformulated as the problem $\inf_{x\in \mathcal D\subseteq \left\{-1,~1\right\}^n}(c^Tx+\gamma)/(d^Tx+\mu)$ with $\gamma\in\Bbb R$, $\mu\in\Bbb R$, though there is a minor difference in the expression. In \cite{R1998,R2013}, Radzik applied the Newton-Dinkelbach method to show that {\rm (LFC)} can be solved in $O(n^2\log^2n)$. Then in \cite{WYZ2006}, Wang et. al. improved the result to $O(n^2\log n)$. Recently in \cite{DZBL2023}, Dadush et. al. presented an accelerated version of the Newton-Dinkelbach method and showed that the iteration bound could be improved to $O(n\log n)$. 
\subsection{Polynomially solvable cases for unconstrained ($-$1,1)-quadratic programming problem}
 In {\rm (P)}, when $B=0,~\beta=1$, {\rm (P)} reduces to an unconstrained ($-$1,1)-quadratic programming problem:
\begin{eqnarray*}
{\rm (QP)}~&\min\limits_{}& x^TQx \\
&{\rm s.t.}& x\in \left\{-1,~1\right\}^n,
\end{eqnarray*}
where $Q$ is a rational $n\times n$ symmetric matrix. {\rm (QP)} contains the max-cut problem, which is NP-hard, as a special case \cite{GJ1979}. It further implies that {\rm (P)} is NP-hard. {\rm (QP)} has many applications in financial analysis \cite{MY1980}, cellular radio channel assignment \cite{CS1995}, statistical physics and circuit design \cite{BGJR1988,GJR1989,P1984}. 

Many polynomially solvable cases of {\rm (QP)} have been identified. Allemand et. al. showed that when $Q\preceq 0$ is a negative semidefinite matrix with fixed rank $p$ and its spectral decomposition is explicitly given, then with the help of zonotope, {\rm (QP)} can be solved in $O(n^{p-1})$ for $p\geq 3$ and $O(n^p)$ for $p\leq 2$ \cite{AFLS2001}. Then it is extended to the condition which is related to the diagonal entries of the matrix $Q$, which is stated in the following two propositions:
\begin{proposition}[\cite{BN2011}]\label{prop1}
For fixed integers $p\geq 2$, if the matrix $Q$ (given by its nonzero eigenvalues and associated eigenvectors) has rank at most $p$ and nonpositive diagonal entries only, then problem {\rm (QP)} can be solved in time $O(n^{p-1}\log(n))$.
\end{proposition}
\begin{proposition}[\cite{BN2011}]\label{prop2}
For a fixed integers $p\geq 2$, if the matrix $Q$ (given by its nonzero eigenvalues and associated eigenvectors) has rank at most $p$ and $O(\log(n))$ positive diagonal entries, then problem {\rm (QP)} can be solved in $O(n^{p}\log(n))$.
\end{proposition}

\begin{proposition}\label{rem1}
If $rank(Q)=1$, and the matrix $Q$ (given by its nonzero eigenvalues and associated eigenvectors) has $O(\log(n))$ positive diagonal entries or has nonpositive diagonal entries only, then problem {\rm (QP)} can be solved in  $O(n)$.\footnote{
In \cite{AFLS2001}, the authors claim that `` when $rank(Q)=1$, {\rm (QP)} is polynomially solvable''. However, this is not always true. When $rank(Q)=1,~Q\preceq 0$, {\rm (QP)} is NP-hard since partition problem, which asks whether the linear equation $\{a^Tx=0,~x\in\{-1,1\}^n\}$ has a solution for any given integer vector a, is NP-hard.
The NP-complete partition problem can be answered by solving a specific instance of {\rm (QP)}, wherein $Q = aa^T$.
}
\end{proposition}

 The techniques for solving {\rm (LFC)} and the polynomially solvable cases of {\rm (QP)} inspire us to give a construction of polynomial complexity of {\rm (P)}.
The paper is organized as below. In Section 2, we equivalently reformulate {\rm (P)} as the problem of finding the unique root of a parametric function and solve it with the accelerated Newton-Dinkelbach method. In each subproblem, we need to solve an unconstrained ($-1$,1)-quadratic programming problem. Finally, we derive a polynomially solvable case of {\rm (P)}. Conclusions are made in Section 3.

%
  
\section{Polynomially solvable case of {\rm (P)}}
\subsection{\textbf{The classical Newton-Dinkelbach method}}
It is well-known that {\rm (P)} is equivalent to finding the unique root of $f(\delta)$:
\begin{eqnarray}\label{ND}
f(\delta)=\min_{x\in\left\{-1,1\right\}^n} \left\{x^TAx+\alpha -\delta(x^TBx+\beta)\right\},
\end{eqnarray}
which is a continuous, concave, strictly decreasing, piecewise linear function \cite{DZBL2023,R1998,R2013}. As we have assumed that $\delta^*\geq 0$, which is the finite optimal objective value of {\rm (P)}. It is equivalent to $f(\delta^*)=0$. 

The classical Newton-Dinkelbach algorithm for solving $f(\delta)$ starts with a starting point $\delta^1$ and a supergradient $g^1\in\partial f(\delta^1):=\left\{g:f(\delta)\leq f(\delta_1)+g(\delta-\delta_1),~\forall \delta\in\Bbb R\right\}$ such that $\delta_1\geq0$, $f(\delta^1)\leq 0$ and $g^1<0$. In each iteration $i\geq 1$, the algorithm maintains a point $\delta^i$, a supergradient $g^i\in\partial f(\delta^i)$, and the function value $f(\delta^i)$. Then we can update $\delta^{i+1}=\delta^i-f(\delta^i)/g^i$.
The algorithm terminates when $f(\delta^i)=0$. 

\begin{figure}[htbp]
\centering
\includegraphics[width=0.6\textwidth]{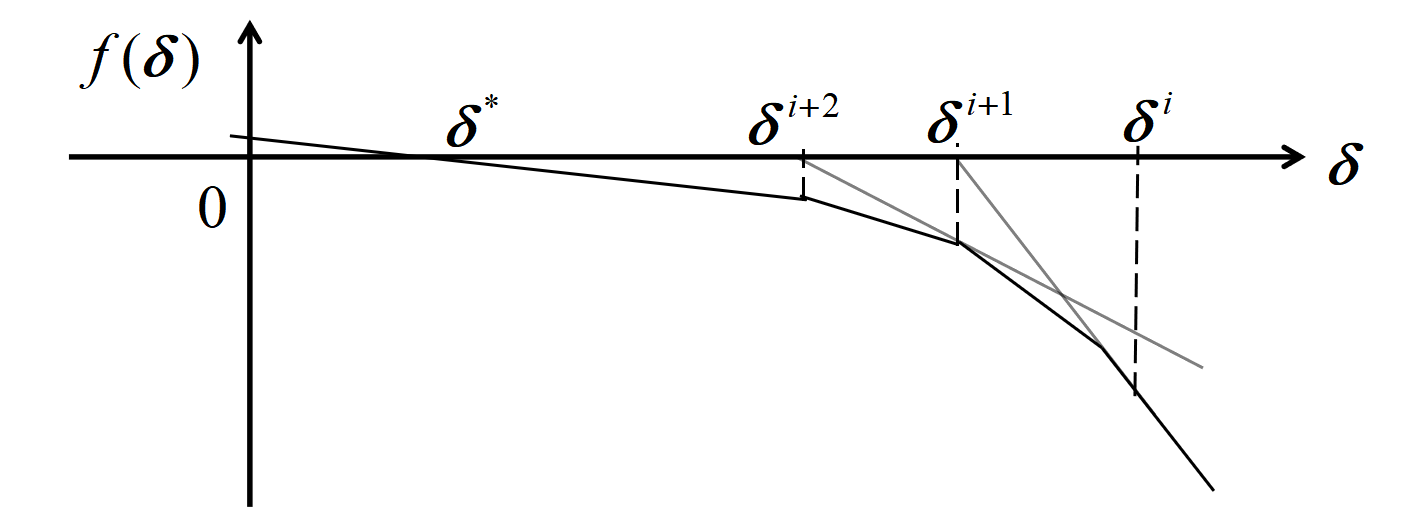}
\caption{The Newton-Dinkelbach method for solving $f(\delta)=0$.}
\label{fig1}
\end{figure}

%
%
It can be proved in the following Lemma \ref{lem1} that 
 $\delta^i$ is monotonically decreasing, $f(\delta^i)$ is monotonically increasing. The lemma also states that $\|f(\delta^i)\|$ or $\|g^i\|$ decreases geometrically. We can also refer to Figure \ref{fig1} as an example.

\begin{lemma}[\cite{DZBL2023}]\label{lem1}
For every iteration $i\geq 2$, we have $\delta^*\leq\delta^i<\delta^{i-1}$, $f(\delta^*)\geq f(\delta^i)>f(\delta^{i-1})$, and $g^i\geq g^{i-1}$, where the last inequality holds at equality if and only if $g^i=\inf_{g\in\partial f(\delta^i)} g$, $g^{i-1}=\sup_{g\in\partial f(\delta^{i-1})}g$, and $f(\delta^i)=0$. Moreover,
\begin{eqnarray}\label{ineq2}
\frac{f(\delta^i)}{f(\delta^{i-1})}+\frac{g^i}{g^{i-1}}\leq 1.
\end{eqnarray}
\end{lemma}

For further analysis, we apply the Bregman divergence associate with $f(\delta)$, which is defined as: 
\begin{definition}[\cite{DZBL2023}]
  Given a proper concave function $f:\Bbb R\rightarrow \bar{\Bbb R}$, the Bregman divergence associated with $f$ is defined as 
  \begin{equation*}
  D_f(\delta',\delta):=
  \begin{cases}
  f(\delta)+\sup_{g\in\partial f(\delta)}g(\delta'-\delta)-f(\delta')~~~if~\delta\neq \delta',\\
  0~~~~~~~~~~~~~~~~~~~~~~~~~~~~~~~~~~~~~~~~~~~otherwise,
  \end{cases}
  \end{equation*}
  for all $\delta,~\delta'\in dom (f)$ such that $\partial f(\delta)\neq \emptyset$.
\end{definition}

We can see that the Bregman divergence is nonnegative since $f$ is concave. If we apply the Bregman divergence associate with $f(\delta)$ to analyze the classical Newton-Dinkelbach method, we can see from the following lemma that $D_f(\delta^*,\delta^{i})$ is monotonically decreasing except in the final iteration, where it may remain unchanged.

\begin{lemma}[\cite{DZBL2023}]\label{nonnegative1}
For every iteration $i\geq 2$, we have $D_f(\delta^*,\delta^i)\leq D_f(\delta^*,\delta^{i-1})$, which holds at equality if and only if
$
g^{i-1}=\sup_{g\in\partial f(\delta^{i-1})}g 
$
and 
$f(\delta^i)=0$.
\end{lemma}
%
%

\subsection{\textbf{An accelerated Newton-Dinkelbach method and its computational complexity}}
In order to accelerate the Newton-Dinkelbach method for solving (\ref{ND}), 
we follow the method proposed for solving {\rm (LFC)}  by Dadush et. al \cite{DZBL2023} and perform an aggressive guess that $\delta^\prime:=2\delta-\delta^i$ on the next point. The corresponding detail is shown in lines 7-9 of Algorithm 1. We follow the name from \cite{DZBL2023} and call the algorithm ``Look-Ahead Newton''. If the ``look-ahead'' guess is successful, then the algorithm is accelerated significantly, otherwise, we are not too far away from the optimal solution.

\begin{algorithm}[H]\label{alg1}
    \caption{Look-ahead Newton method for finding the root of $f(\delta)$. 
    }
    \label{alg2}
    \textbf{Input:} An initial point $\delta^{(1)}\geq 0$, supergradient $g^{(1)}\in\partial f(\delta^{(1)})$, where $f(\delta^{(1)})\leq 0$ and $g^{(1)}<0$.\\
    \textbf{Output:} The optimal solution $\delta^*$.
    \begin{algorithmic}[1]
        \State $i\leftarrow 1$
        \While{$f(\delta^{(i)})<0$}
        \State$\delta\leftarrow\delta^{(i)}/g^{i}$
        \State$g\in\partial f(\delta)$
        \State $\delta^\prime\leftarrow 2\delta-\delta^{(i)}$
        \State $g^\prime\in\partial f(\delta^\prime)$
        \If{$f(\delta^\prime)<0$ and $g^\prime<0$}
        \State $\delta\leftarrow\delta^{\prime}$, $g\leftarrow g^\prime$
        \EndIf
        \State $\delta^{(i+1)}\leftarrow \delta$, $g^{(i+1)}\leftarrow g$.
        \State $i\leftarrow i+1$
        \EndWhile
        \State \textbf{return} $\delta^{(i)}$
    \end{algorithmic}        
\end{algorithm}


We follow the technology line proposed by Daudush et. al. \cite{DZBL2023} to employ the Bregman divergence associated with $f(\delta)$ to analyze the computational complexity of Algorithm 1.

 The next lemma demonstrates the advantage of using the look-ahead Newon method. 

\begin{lemma}[\cite{DZBL2023}]\label{dec1}
For every iteration $i>2$ in Algorithm 1, we have $D_f(\delta^*,\delta^i)<\frac{1}{2}D_f(\delta^*,\delta^{i-2})$.
\end{lemma}

\begin{remark}
It is easy to verify that Lemmas \ref{lem1}, \ref{nonnegative1} also hold for Algorithm 1. Besides, Lemmas \ref{lem1}, \ref{nonnegative1} and \ref{dec1} are all proposed by Daudush et. al. in \cite{DZBL2023} for finding the root of the following parametric function: 
\[\tilde{f}(\tilde{\delta})=\inf\left\{(c-\tilde{\delta}d)^Tx:x\in\mathcal D\right\}\]
corresponding to {\rm (LFC)} with Dinkelbach-Newton method or look-ahead Newton method. It is easy to verify that these three lemmas also hold for $f(\delta)$ with respect to {\rm (P)}.
\end{remark}

In order to analyze the computational complexity of the look-ahead Newton method for finding the root of $f(\delta)$, i.e. (\ref{ND}), we make a transformation of the primal problem (P). Since $x\in\left\{-1,1\right\}^n$, we know that $x^Tx=n$ holds. Then 
\[
x^TAx+\alpha=x^TAx+\frac{x^Tx}{n}\alpha=x^T(A+\frac{1}{n}I_n)x:=x^T\tilde{A}x,
\] 
where $I_n$ is an $n\times n$ identity matrix. Similarly, $x^TBx+\beta=x^T(B+\frac{1}{n}I_n)x:=x^T\tilde{B}x$.
So we can obtain the following equivalent problem of (P):
\begin{eqnarray}
(PE1)~~&\min& \frac{x^T\tilde{A}x}{x^T\tilde{B}x} \\
&{\rm s.t.}& x\in \left\{-1,~1\right\}^{n},
\end{eqnarray}
and the corresponding parametric function is
\begin{eqnarray}\label{ND1}
f(\delta)=\min_{x\in\left\{-1,1\right\}^n} \left\{x^T\tilde{A}x-\delta x^T\tilde{B}x\right\},
\end{eqnarray}
which is equal to (\ref{ND}). 

In \cite{R2013}, Radzik proposed a bound on the length of a geometrically decreasing sequence of sums of numbers: 
\begin{lemma}[\cite{R2013}]\label{lin1}
Let $c=(c_1,c_2,...,c_p)$ be a $p$-dimensional vector with nonnegative real components, and let $y_1,y_2,...,y_q$ be vectors from $\left\{-1,0,1\right\}^p$. If for all $i=1,2,...,q-1$, $0<y_{i+1}^Tc\leq \frac{1}{2}y_i^Tc$, then $q=O(p\log p)$.
\end{lemma}

We make an extension of Lemma \ref{lin1} and obtain the following result, which is the main tool in analyzing the computational complexity of the look-ahead Newton method for finding the unique root of $f(\delta)$, i.e. (\ref{ND1}).

\begin{lemma}\label{lin2}
Let $C=(c_{ij})_{n\times n}$ be an $n\times n$-dimensional symmetric matrix. Let $y_1,y_2,...,y_q$ be vectors from $\left\{-1,0,1\right\}^{n^2}$. If for all $i=1,2,...,q-1$, $0<y_{i+1}^TCy_{i+1}\leq \frac{1}{2}y_i^TCy_i$ holds, then $q=O(n^2\log n)$.
\end{lemma}
{\bf {Proof.}}
Since $0<y_{i+1}^TCy_{i+1}\leq \frac{1}{2}y_i^TCy_i$ holds for $\forall~i=1,2,...,q-1$, then we can take the trace operation and obtain
\begin{eqnarray}
&&0<tr(y_{i+1}^TCy_{i+1})\leq \frac{1}{2}tr(y_i^TCy_i)~~~~~~~~~~~~~~~~~~~~~~~~~~~~~~~~\nonumber\\
&\Rightarrow&0<tr(Cy_{i+1}y_{i+1}^T)\leq \frac{1}{2}tr(Cy_iy_i^T)\nonumber\\
(by~letting~Y_i=y_iy_i^T)~&\Rightarrow&0<tr(CY_{i+1})\leq \frac{1}{2}tr(CY_i)\nonumber\\
(by~stacking~the~columns~of~C~and~Y_{i})~&\Rightarrow&0<\tilde{c}^T\tilde{y}_{i+1}\leq\frac{1}{2}\tilde{c}^T\tilde{y}_{i}\label{eq2}\\
(by~letting~ \tilde{c}^T\tilde{y}_{i}=|\tilde{c}|^T\tilde{z}_{i})~&\Rightarrow&0<|\tilde{c}| ^T\tilde{z}_{i+1}\leq\frac{1}{2}|\tilde{c}|^T\tilde{z}_{i},\label{eq3}
\end{eqnarray}
where (\ref{eq2}) holds by stacking the columns of $C$, $Y_{i}$ and denoted by $\tilde{c}\in\Bbb R^{n^2}$, $\tilde{y}_i\in\Bbb R^{n^2}$, respectively. Consequently, in (\ref{eq3}), we know that $|\tilde{c}|\geq 0$, $\tilde{z}_i\in\left\{-1,0,1\right\}^{n^2}$. Then according to Lemma \ref{lin1}, we can conclude that $q=O(n^2\log n)$.
{\hfill \quad{$\Box$} }\\

\begin{theorem}\label{thm1}
Algorithm 1 converges in $O(n^2\log n)$ iterations for finding the unique root of $f(\delta)$, i.e. (\ref{ND1}).
\end{theorem}
{\bf {Proof.}}
We know that Algorithm 1 terminates in a finite number of iterations since $f(\delta)$ is piecewise linear. Let $\delta^1>\delta^2>...>\delta^k=\delta^*$ denote the sequence of iterates at the start of Algorithm 1. Because $f$ is concave, we have $D_f(\delta^*,\delta^i)\geq 0$ for all $i=1,...,k$. 
For each $i$, 
\[
\partial f(\delta^i)=\left\{-x^{i^T}\tilde{B}x^i~\mid~x^i\in\arg\min_{x\in \left\{-1,1\right\}^n} (x^T\tilde{A}x-\delta^i x^T\tilde{B}x)\right\}.
\]
As $f(\delta^*)=0$, the Bregman divergence of $\delta^i$ and $\delta^*$ can be written as 
\begin{eqnarray}
D_f(\delta^*,\delta^i)&=&f(\delta^i)+\max_{g\in\partial f(\delta^i)}g(\delta^*-\delta^i)\nonumber\\
&=&x^{i^T}(\tilde{A}-\delta^i\tilde{B})x^i-x^{i^T}\tilde{B}x^i(\delta^*-\delta^i)\label{eq21}\\
&=&x^{i^T}(\tilde{A}-\delta^*\tilde{B})x^i,\nonumber
\end{eqnarray}
where in (\ref{eq21}), we pick $x^i=\tilde{x}^i$ such that $\max_{g\in\partial f(\delta^i)}g=-\tilde{x}^{i^T}\tilde{B}\tilde{x}^i$. Then according to Lemma \ref{dec1}, 
\[
x^{i^T}(\tilde{A}-\delta^*\tilde{B})x^i=D_f(\delta^*,\delta^i)<\frac{1}{2}D_f(\delta^*,\delta^{i-2})
=\frac{1}{2}(x^{i-2})^T(\tilde{A}-\delta^*\tilde{B})x^{i-2},
\]
where $x^i\in\left\{-1,1\right\}^n\subseteq \left\{-1,0,1\right\}^n$ holds for all $3\leq i\leq k$.  By Lemma \ref{nonnegative1}, we also know that $D_f(\delta^*,\delta^i)>0$ for all $1\leq i\leq k-2$. 
 Then it follows from Lemma \ref{lin2} that $k=O(n^2\log n)$. The proof is complete.
{\hfill \quad{$\Box$} }\\

\subsection{\textbf{Polynomially solvable cases for solving the ($-$1,1)-quadratic optimization problem}}
In each iteration, for a given $\bar{\delta}$, it must hold that $\bar{\delta}
\geq \delta^*\geq 0$. And we need to solve the following ($-$1,1)-quadratic programming problem:
\begin{eqnarray*}\label{BQ1}
\min_{x\in\left\{-1,1\right\}^n} \left\{x^T(A-\bar{\delta}B)x+\alpha-\bar{\delta}\beta\right\},
\end{eqnarray*}
which is equivalent to
\begin{eqnarray}\label{BQ2}
{\rm (QP(\bar{\delta}))}~~\min_{x\in\left\{-1,1\right\}^n} \left\{x^T(A-\bar{\delta}B)x\right\}
\end{eqnarray}
in the sense that they share the same optimal solution. Ben-Ameur and Neto  \cite{BN2011} have demonstrated the condition under which the unconstrained ($-$1,1)-quadratic programming problem {\rm (QP)} can be solved in polynomial time in Propositions \ref{prop1} and \ref{prop1}. Based on which we can obtain the condition for solving {\rm (QP($\bar{\delta}$))} in polynomial time.


\begin{theorem}\label{thm2}
For fixed integers $r_a$ and $r_b$, if the matrices $A$ and $B$ (given by their nonzero eigenvalues and associated eigenvectors, respectively) have rank at most $r_a$ and $r_b$, respectively, $A$ has nonpositive diagonal entries only, $B$ has nonnegative diagonal entries only, then problem {\rm (QP($\bar{\delta}$))} can be solved in time $O(n^{r_a+r_b-1}\log(n))$.
\end{theorem}
{\bf {Proof.}}
We denote the diagonal entries of matrix $A\in\Bbb R^{n\times n}$ by $A_{ii},~i=1,...,n$ and the same symbol is also used for other matrices. Since $rank (A)\leq r_a$, $rank (B)\leq r_b$. It holds that $rank(A-\bar{\delta}B)\leq r_a+r_b$. Since $A_{ii}\leq 0,~i=1,2,...,n$, $B_{ii}\geq 0,~i=1,2,...,n$, $\bar{\delta}\geq 0$. It holds that $(A-\bar{\delta}B)_{ii}\leq 0,~i=1,2,...,n$. Subsequently, according to Proposition \ref{prop1} (or Proposition \ref{rem1} when $rank(A-\bar{\delta}B)=1$), {\rm (QP($\bar{\delta}$))} can be solved in time $O(n^{r_a+r_b-1}\log(n))$. 
{\hfill \quad{$\Box$} }\\
%

\begin{theorem}\label{thm3}
For fixed integers $r_a$ and $r_b$, if the matrices $A$ and $B$ (given by their nonzero eigenvalues and associated eigenvectors, respectively) have rank at most $r_a$ and $r_b$, respectively, $A$ has $O(\log(n))$ positive diagonal entries, $B$ has $O(\log(n))$ negative diagonal entries, then problem {\rm (QP($\bar{\delta}$))} can be solved in $O(n^{r_a+r_b}\log(n))$.
\end{theorem}
{\bf {Proof.}}
We use $\#\left\{A_{ii}>0,~i=1,2,...,n\right\}$ to denote the number of positive entries on the diagonal of the matrix $A\in\Bbb R^{n\times n}$, and the same symbol is also used for other matrices. Since $rank (A)\leq r_a$, $rank (B)\leq r_b$. Then $rank(A-\bar{\delta}B)\leq r_a+r_b$. Since $\#\left\{A_{ii}>0,~i=1,2,...,n\right\}\leq O(\log(n))$, 
$\#\left\{B_{ii}<0,~i=1,2,...,n\right\}\leq O(\log(n))$,
$\bar{\delta}\geq 0$. It holds that $\#\left\{(A-\bar{\delta}B)_{ii}>0,~i=1,2,...,n\right\}\leq O(\log(n))$.
Then according to Proposition \ref{prop2} (or Proposition \ref{rem1} when $rank(A-\bar{\delta}B)=1$), {\rm (QP($\bar{\delta}$))} can be solved in $O(n^{r_a+r_b}\log(n))$. 
{\hfill \quad{$\Box$} }\\

\subsection{\textbf{Polynomially solvable cases for solving {\rm (P)}}}
In order to solve {\rm (P)}, we separate the procedure into two steps. (i) We employ the look-ahead Newton method to solve the parametric function and the computational complexity is $O(n^2\log n)$, see Theorem \ref{thm1}. (ii) In each iteration of the look-ahead Newton algorithm, we need to solve a ($-$1,1)-quadratic programming problem. According to Theorems \ref{thm2} and \ref{thm3}, we obtain the cases under which the ($-$1,1)-quadratic programming problem can be solved in $O(n^{r_a+r_b-1}\log n)$ ($r_a$ and $r_b$ are fixed integers) or $O(n^{r_a+r_b}\log n)$, respectively.

%

Based on the derivations above, we can directly obtain polynomially solvable cases of unconstrained ($-$1,1)-quadratic fractional optimization, i.e. {\rm (P)}:  

\begin{theorem}
The problem {\rm (P)} can be solved in $O(n^{r_a+r_b+1}\log^2 n)$ if for fixed integers $r_a$ and $r_b$, the matrices $A$ and $B$ (given by their nonzero eigenvalues and associated eigenvectors, respectively) have rank at most $r_a$ and $r_b$, respectively, and $A$ has nonpositive diagonal entries only, $B$ has nonnegative diagonal entries only.
\end{theorem}
\begin{theorem}
The problem {\rm (P)} can be solved in $O(n^{r_a+r_b+2}\log^2 n)$ if for fixed integers $r_a$ and $r_b$, the matrices $A$ and $B$ (given by their nonzero eigenvalues and associated eigenvectors, respectively) have rank at most $r_a$ and $r_b$, respectively, and $A$ has $O(\log(n))$ positive diagonal entries, $B$ has $O(\log(n))$ negative diagonal entries.
\end{theorem}

\section{Conlusion}
In this paper, in order to solve an unconstrained ($-$1,1)-quadratic fractional optimization in the form as shown in {\rm (P)}, we equivalently reformulate it as finding the root of a parametric function and apply the look-ahead Newton method to solve it. In each iteration, an unconstrained ($-$1, 1)-quadratic optimization subproblem needs to be solved. 
We show that the problem {\rm (P)} (where $A$ and $B$, given by their nonzero eigenvalues and associated eigenvectors, have ranks not exceeding fixed integers $r_a$ and $r_b$, respectively.)
 can be solved in $O(n^{r_a+r_b+1}\log^2 n)$ by the accelerated Newton-Dinkelbach method when the matrices $A$ has nonpositive diagonal entries only, $B$ has nonnegative diagonal entries only. Furthermore, the problem {\rm (P)} can be solved in $O(n^{r_a+r_b+2}\log^2 n)$ when $A$ has $O(\log(n))$ positive diagonal entries, $B$ has $O(\log(n))$ negative diagonal entries.
In the future, we will further consider an unconstrained ($-$1, 1) fractional programming problem when the quadratic function in both the numerator and the denominator includes linear terms.
\section*{Acknowledgments}
This research was supported by the National Natural Science Foundation of China under grants 12101041, 12171021.

\end{document}